\documentstyle[a4wide,12pt]{article}
\def\sym{\fam\comfam\com}
\font\tensym=msbm10 \font\sevensym=msbm7 \font\fivesym=msbm5
\newfam\symfam
\textfont\symfam=\tensym \scriptfont\symfam=\sevensym
\scriptscriptfont\symfam=\fivesym
\def\sym{\fam\symfam\relax}
\def\N{{\sym N}}
\def\Z{{\sym Z}}

\def\R{{\sym R}}

\def\X{{\sym X}}
\def\E{{\sym E}}
\def\P{{\sym P}}
\def\1{{\bf 1}}
\newcommand{\nit}{\noindent}
\newcommand{\ti}{\partial}
\newcommand{\bh}{\bar H}
\newcommand{\ab}{\bg<}
\newcommand{\ba}{\bg>}

\newcommand{\cl}{{\cal C}}
\newcommand{\dl}{{\cal D}}
\newcommand{\mm}{{\cal M}_+}
\newcommand{\jj}{{\cal J}_H}
\newcommand{\ii}{{\cal I}}
\newcommand{\tf }{{\cal T}}
\newcommand{\Cf}{{\cal B}}               \newcommand{\fkn}{\f{k}{N}}
\newcommand{\bl}{\label}
\newcommand{\fxn}{\bg(\f{x}{N}\bg)}
\newcommand{\Bq}{\begin{equation}}
\newcommand{\Eq}{\end{equation}}
\newcommand{\bg}{\bigg}
\newcommand{\fl}{\frac{1}{2l+1}}
\newcommand{\fN}{\f{1}{N}}
\newcommand{\Bg}{\Bigg}
\newcommand{\pro}{\P_{\rho,p}^N}
\newcommand{\hpro}{\P_{\ma,N}^{p,H}}
\newcommand{\esp}{\E^N_{\rho,p}}
\newcommand{\f}{\frac}
\newcommand{\lspA}{\overline{\lim_{A\ra\ty}}\ }
\newcommand{\lspN}{\overline{\lim_{N\ra\ty}}\ }
\newcommand{\lifN}{\liminf_{N\ra\ty}\ }
\newcommand{\lspl}{\overline{\lim_{l\ra\ty}}\ }
\newcommand{\lspk}{\overline{\lim_{k\ra\ty}}\ }

\newcommand{\lspdt}{\overline{\lim_{\dt\ra 0}}\ }
\newcommand{\nbr}{{\bar\nu}_\rho}
\newcommand{\nbrl}{{\bar\nu_{\rho,0,l}^p}}
\newcommand{\nrl}{{\bar\nu_{\rho,0,l}}}
\newcommand{\nbrxl}{{\bar\nu_{\rho,x,l}}}
\newcommand{\lspe}{\overline{\lim_{\e\ra 0}}\ }
\newcommand{\lspfa}{\overline{\lim_{\fa\ra 0}}\ }
\newcommand{\Ww}{W_{N,\e}^{H,\Psi}}
\newcommand{\sx }{\sum_x}
\newcommand{\wg}{\Lambda}
\newcommand{\ld}{\lambda}                     \newcommand{\ta}{\theta}
\newcommand{\fa}{\alpha}                      \newcommand{\ma}{\gamma}
\newcommand{\lm}{\lambda}
\newcommand{\e}{\varepsilon}
               \newcommand{\bt}{\beta}
\newcommand{\vp}{\varphi}                     \newcommand{\dt}{\delta}

\newcommand{\eN}{(2\e N+1)}
\newcommand{\om}{\omega}                \newcommand{\ty}{\infty}
\newcommand{\ra}{\rightarrow}

           \newcommand{\di} {{\ \rm
d}}                 

               \newcommand{\fx}{\vp
p_x^{-1}}                  
\newcommand{\fk}{f_{x,y,l}^p}

                           \newcommand{\si}{\sigma}
\newcommand{\pn}{\pi_t^N}
\newcommand{\feN}{\frac{1}{\eN}}
\newcommand{\A}{A^{l,k,\dt}}           
       \newcommand{\vk}{\vskip0.3truecm}

\begin{document}
\title{\large{Large deviations for a zero mean asymmetric
zero range process in random media.}}\vskip0.2cm
\author{A. Koukkous\kern -0.5pt
\addtocounter{footnote}{10}
\renewcommand{\thefootnote}{\alph{footnote}}\footnotemark%
\ \ and\ H. Guiol\kern -2pt \addtocounter{footnote}{-5
}\renewcommand{\thefootnote}{\alph{footnote}} \footnotemark}
\maketitle
\renewcommand{\thefootnote}{\alph{footnote}}
\addtocounter{footnote}{7}\footnotetext{$^{,k}$IMECC-UNICAMP, P.B.
6065, 13083-970, Campinas, SP, Brasil.}

\begin{abstract}
We consider an asymmetric zero range process in infinite volume with zero mean 
 and random jump rates starting from equilibrium.
We investigate the large deviations from the hydrodynamical limit of the
empirical distribution of particles and prove an upper and a lower
bound for the large deviation principle. Our main argument is based
on a super-exponential estimate in infinite volume. For this we extend to our case a method developed by Kipnis \& al. (1989) and Benois \& al. (1995).
\vk\nit {\it Keywords: } Asymmetric zero range,
Hydrodynamical limit, Large deviations, Random environment. \vk\nit {\it
AMS 2000 classification } 60K35, 60K37, 82C22.
\end{abstract}
\noindent


\section{Introduction}
The so called zero range process is one of the simplest particle systems that
has been systematically and successfully investigated in random or
inhomogeneous media in the last few years (see for instance Benjamini
\& al. (1996), Evans (1996),  Krug-Ferrari (1996), Landim (1996),
Gielis \& al. (1998), Bahadoran (1998), Sepp\"{a}l\"{a}inen-Krug
(1999), Koukkous (1999), Andjel \& al. (2000)) . 

The zero process can be described informally as follows.
Particles are distributed on the $d$-dimensional lattice $\Z^d$.
Each particle at site $x$ of $\Z^d$ jumps, with a rate depending only
on the total number of particles standing at this site, to the left or
to the right.  

In what follows, 
we consider a sequence of random variables $p=(p_x)_{x\in\Z^d}$ ({\it called an
environment}) in $[a_0,a_1]$ (where $0<a_0\leq a_1<\ty$). According to
$p$ the jump rates of the process are accelerated or decelerated by
the value $p_x$ at site $x$.

Benjamini \& al. (1996) have studied the asymmetric version of a
zero-range process in infinite volume when the environment is an
i.i.d. sequence of random variables (with $a_1=1$) and have proved
the asymptotic hydrodynamical behavior of the system. 
Koukkous (1999) proved the hydrodynamical limit in the symmetric
case for a stationary and ergodic environment whose marginal law
is absolutely continuous with respect to the Lebesgue measure. 
In particular he showed that the empirical measure of particles converges in
probability to the weak solution of a non-linear diffusion
equation which does not depend on the environment $p$ and
generalized in this way some results of Benjamini \& al. (1996).

The equilibrium fluctuations (Central limit results for the density
field) were studied in G. Gielis \& al. (1998). They proved that the
density field converges weakly to a generalized Ornstein-Uhlenbeck
process.

Recently, Andjel \& al. (2000) showed the convergence to the maximal
invariant measure for an asymmetric zero range process with
constant rate in inhomogeneous or random media in dimension 1
starting from an upper-critical non-equilibrium measure.

\vk In this spirit of hydrodynamical behavior investigation, a
natural open question can be formulated as follows: From the
hydrodynamical limit of the empirical measure with some continuous
density $\mu(\cdot)$ (with respect to Lebesgue measure) and given an
event $\Gamma$ for which $\mu\notin\bar\Gamma$, how to control the
``{\it deviant}" behavior of the system inside  $\Gamma$ ? This
is the subject of large deviation principles (LDP) related to
hydrodynamical limit of the empirical measure.\vk

In this paper, we investigate this question for a $d$-dimensional zero mean
asymmetric zero-range process in random media. In the
deterministic case, the LDP results have been treated by many 
authors among which Landim (1992), Benois (1996) and Benois \&
al. (1995). In this last article Benois \& al. gives an upper and a
lower bound of the LDP in infinite volume for the empirical density
when the process starts from equilibrium. The crucial ingredient of
their arguments focuses on the so-called super-exponential estimate:
it consists to approximate, by some rigorous functions of the density
field, the correlation field obtained by computing some exponential
martingales related to the jumps of particles (see Kipnis \&
al. (1989) and Donsker-Varadhan (1989)). Once one prove this result,
the LDP (and also the hydrodynamical limit) for the empirical measure
is obtained by standard arguments. 

\vk In random environment, the difficulty relies on the absence of
translation invariance of the invariant measures of the process. 
For this reason our approach will also use some results of
Koukkous (1999). 

\vk The paper is organized as follows: in Section 2 we introduce the notation
and assumptions used along the paper and state the main results. 
Section 3 is devoted to the proof of the super-exponential estimate. 
In the last section we give a proof of the upper bound of LDP result. 
We omit the proof of lower bound since, once one has proven the upper
bound, it is similar to the arguments given in Benois \& al. (1995)
without major modifications.
 \vk
\section{Notation and results}
 \bl{NR}
Let $0<a_0\leq a_1<\infty$ and consider a sequence of random
variables $\{p_x,\ x\in\Z^d\}$ on $[a_0,a_1]$ distributed
according to an ergodic stationary measure $m$, such that its
one-dimensional marginal law is absolutely continuous with respect
to the Lebesgue measure. We assume that $m\{p: a_0\leq p_0\leq
a_1\}=1$ and for every $\e>0$, $m\{p: p_0\in [a_0,a_0+\e)\}m\{p:
p_0\in (a_1-\e,a_1]\}>0$.

\vk We denote by $\X_d:=\N^{\Z^d}$ the configuration space and by
Greek letters $\eta$ and $\xi$ its elements. As usual $\eta(x)$
stands for the total number of particles at site $x$ for
the configuration $\eta$. For each environment $p$, we are interested
in the Markov process ${(\eta_t)}_{t\geq 0}$ on $\X_d$ whose
generator is defined by
\Bq \label{generateur}
({\cal
L}_pf)(\eta)=\sum_{\scriptstyle{x,y\in\Z^d}}
p_xg(\eta(x))T(x,y)[f(\eta^{x,y})-f(\eta)],
\Eq
where
$f:\X_d\rightarrow\R$ is a bounded cylinder function,
that is $f$ only depends on $\eta$ through a finite number of coordinates.
$T(\cdot,\cdot)$ is a transition probability on $\Z^d$. The
function $g$ is positive and vanishes at $0$: 
$g(0)=0 < g(k)\ \mbox{for all}\ k\geq 1$. In the previous formula,
$\eta^{x,y}(z)$ is the configuration obtained from $\eta$ when a
particle jumps from $x$ to $y$:
\[
\eta^{x,y}(z) =
\left\{\begin{array}{ll}
                       \eta(z) & \ \mbox{if $z\neq x,y$}  \\
                       \eta(x)-1& \ \mbox{if $z=x$}  \\
                       \eta(y)+1& \ \mbox{if $z=y$\ \ \ .}
                         \end{array}\right.
\]
For every non-negative real $\vp$ we denote by $\nu_\vp^p$ the
product measure on $\X_d$ whose marginals are defined by
\[
\nu_\vp^p\{\eta : \eta(x) =k\} = \frac{1}{Z(\fx)}
\frac{(\fx)^k}{g(k)!}, \ \ \ \mbox{for all}\ \  k\geq 0,
\]
where $g(k)!=g(1)g(2)...g(k)$ if $k>0$ and $g(0)!=1$. Under some
hypotheses (for instance  [H1] and [H2] in what follows), those
measures  are invariant for the 
process. In this formula, $Z : \R_+\ra\R_+$ is the partition
function
\[
Z(\vp)=\sum_{k \geq 0} \frac{\vp^k}{g(k)!}.
\]
Let $\vp^*$ be the radius of convergence of $Z(\cdot)$; we assume
that
 \Bq
 \bl{eq1}
\lim_{\vp\uparrow\vp^*}Z(\vp) =+\infty.
 \Eq
 Denote by $\nu_{\vp}(\cdot):=\nu_\vp^1(\cdot)$ the invariant measure of the process
$(\eta_t)_{t\geq 0}$ when $m$ is the Dirac measure concentrated on
the set $\{p : p_x=1,x\in \Z^d\}$ (see Andjel (1982)). We define
$M : [0,\vp^*)\ra\R_+$ by
$M(\vp)=\nu_\vp[\eta(0)]$, the expected number of particles at $0$
with respect to $\nu_{\vp}$.

\vk A simple computation shows that $M(\vp)=\vp\partial_\vp \log
Z(\vp)$ and from assumption (\ref{eq1}) we check that $M$ is an
increasing, continuous, one-to-one function from $[0,\vp^*)$ to
$\R_+$.

We define the ``density" of particles ({\it i.e.} the expected
number of particles at $0$) with respect to the random media by the
continuous and increasing function $R : [0,a_0\vp^*)\ra\R^+$ such
that
\[
R(\vp)=m[M(\vp p_0^{-1})]
\]
and in order to ensure the existence of an invariant measure for
any given value of the density, we assume that
 \Bq
 \bl{eq2}
\lim_{\vp\uparrow a_0\vp^*} R(\vp) =\infty.
 \Eq
Under this assumption the function $R$ is one to one from
$[0,a_0\vp^*)$ to $\R_+$. We denote by $\Phi$ its inverse (which
is also a continuous increasing bijection).

\vk For a density $\rho>0$ we write $$\nbr^p=\nu_{\Phi(\rho)}^p.$$

In the following we state all the hypotheses assumed throughout
this paper. 

\vk {\bf [H1]}\ \ The transition probability $T(\cdot,\cdot)$ on
$\Z^d$ is a zero-mean irreducible translation invariant
probability with finite range. That is
\[
\begin{array}{c}
\hspace{2.1cm}T(x,y)=T(0,y-x)=:T(y-x), \ 
\hspace{2.5cm}\\ \mbox{there exists
a constant }A>0\mbox{ such that } T(x)=0\ \ \mbox{if}\ \ |x|\geq
A\\ \ \mbox{and }\ \ \displaystyle{ \sum_{x\in\Z^d}x~T(x)=0}.
\end{array}
\]

{\bf [H2]}\ \  The rate function $g$ has bounded variation:
\[
\hspace{3.1cm}g^*=\sup_k|g(k+1)-g(k)|<\infty.\hspace{4.5cm}
\]

Under the hypotheses [H1] and [H2] there
exists a unique Markov process with corresponding generator
defined by (\ref{generateur}) for the deterministic case {\it
i.e.} $p\equiv 1$ (see Andjel (1982)). Andjel's proof applies also 
in the case we consider.

Let $(\si_{ij})_{\{1\leq i,j\leq d\}}$ be a symmetric nonnegative
definite matrix defined by the covariance matrix of the transition
probability $T(\cdot)$:
\[
\si_{ij}=\sum_{y\in\Z^d}y_iy_jT(y)\ \ \ \mbox{where}\ \
y=(y_1,\cdots,y_d).
\]

{\bf [H3]}\ In order to avoid the degenerate case of the
hydrodynamic equation, we assume $(\si_{ij})_{\{1\leq i,j\leq d\}}$ to be a positive
definite matrix. That is there exists $\kappa>0$ such that
\[\hspace{3.1cm}\sum_{i,j}\si_{ij}x_ix_j\geq\kappa\sum_i
x_i^2,\ \ \mbox{for all}\ \  x=(x_1,\cdots,x_d)\in \R^d.\hspace{3.5cm}
\]

{\bf [H4]} To ensure some finite exponential moments of $\eta(x)$ under the measures $\nu_\vp^p$ we shall assume that there exists a convex and increasing function $\om
:\R_+\longrightarrow\R_+$ such that\\ (i)  $ \om(0)=0,$ \\ (ii)
$\lim_{x\ra\infty}(\f{\om(x)}{x})=\infty$\ and\\ (iii) for all
density $\vp$ there exists a positive constant $\ta:=\ta(\vp)$ such that
\[
\nu_\vp\bigg[\exp\big\{\ta\om(\eta(0))\big\}\bigg]<\infty.
\]
This last assumption ensures also that $Z(\cdot)$ has infinite radius of
convergence. It holds for exemple if $g(k+1)-g(k)\geq g^*_0$ for
some constant $g^*_0$ and $k$ sufficiently large.\par

We will denote by $\om^*$ the Legendre transform of $\om$ given
by:
 \Bq
 \bl{ome}
\om^*(x)=\sup_{\fa>0}\{\fa x-\om(\fa)\}.
 \Eq
In the next paragraphs, we define the state space of the process
and its topology. Denote by ${\cal C}(\R^d)$ (resp. ${\cal
C}_K(\R^d)$) the space of continuous (resp. with compact
support) functions on $\R^d$ with classic uniform norm. Let $\mm$
denote the space of positive Radon measures on $\R^d$ with the weak
topology induced by ${\cal C}_K(\R^d)$ via
$\big<\pi,H\big>:={\large\int} H\di\pi$ for $H\in{\cal C}_K(\R^d)$
and $\pi\in\mm$.

We fix a positive time parameter $\tf >0$. For each realization of
the environment $p$ and all fixed positive density $\rho$, $\pro$
will denote the probability measure on the path space ${\cal
D}([0,\tf ],\X_d)$ corresponding to the Markov process
$(\eta_t)_{t\in[0,\tf ]}$ with generator $N^2{\cal L}_p$ starting
from the measure $\nbr^p$. By $\esp$ we denote the expectation
under $\pro$.

\vk Let $\pi^N_.$ be the empirical measure defined on ${\cal D}([0,\tf
],\mm)$ by
\[
\pn(du)=\f{1}{N^d}\sum_{x\in\Z^d}\eta_t(x)\dt_{x/N}(du),
\]
 for $0\leq t\leq \tf $. Let $Q^N_{\rho,p}$ denote the measure on the path space
$D([0,{\cal T}],\mm)$ associated to the process $\pi_.^N$ with
generator $N^2{\cal L}_p$  starting from $\nbr^p$ .

\vk To investigate the large deviations of the empirical measure,
we shall consider some small perturbations of the zero range
process as mentionned earlier. For this, we will need the following
notation.\par

Let ${\cal C}^{l,k}_K([0,\tf ]\times\R^d)$ denote the space of compact support functions with $l\in\N$ continuous derivatives in time and
$k\in\N$ continuous derivatives in space. Let ${\cal C}_\rho(\R^d)$
be the set defined by
\[
{\cal C}_\rho(\R^d)={\cal C}(\R^d)\cap\{u:\R^d\ra\R^+;\ \ u(x)=\rho \
\mbox{for}\ \ |x|\ \mbox{sufficiently large}\}.
\]

For a fixed $\ma$ in ${\cal C}_\rho(\R^d)$ and for some smooth
function $H$ in $\cl^{1,2}_K([0,\tf ]\times\R^d)$ we consider the
Markov process generated by
\[
N^2({\cal
L}_{N,t}^{p,H}f)(\eta)=N^2\sum_{\scriptstyle{x,y\in\Z^d}}
p_xg(\eta(x))T(y)e^{\{H(t,\f{x+y}{N})-H(t,\f{x}{N})\}}[f(\eta^{x,x+y})-f(\eta)],
\]
where $f$ is a cylinder function. Let $\bar\nu_{\ma,N}^p$ be the
initial product measure of this process with marginals
\[
\bar\nu_{\ma,N}^p\{\eta, \eta(x)=k\}=\bar\nu_{\ma(x/N)}^p\{\eta,
\eta(x)=k\}
\]
for all $x\in\Z^d$ and $k\in \N$. We therefore denote by $
\P_{\ma,N}^{p,H}$ and $Q^{p,H}_{\ma,N}$ the small perturbations of
$\pro$ and $Q^N_{\rho,p}$ respectively.

For any path $\pi_.\in\dl([0,\tf ],\mm)$, denote by $u_t$ the
Radon-Nikodym derivative of $\pi_t$ with respect to the Lebesgue
measure $\lm$: $u_t:=\f{\di\pi_t}{\di\lm}$. Let $\cal A=\cal
A(\rho)$ be the space path of $\pi\in\dl([0,\tf ],\mm)$ such that
$u_t$ is the solution of the PDE
\[
\mbox{\bf(E)} \hspace{3cm}\left\{\begin {array}{ll}
              \partial_t u
&=(\si/2)\triangle(\Phi(u))-\sum_{i=1}^{d}\partial_{x_i}(\Phi(u)\partial_{x_i}H)\\
              u(0,\cdot)&=\gamma(\cdot)\ \ .\hspace{4.5cm}
      \end{array}
\right.\hspace{2cm}
\]

for some $\gamma\in{\cal C}_\rho(\R^d)$ and
some $H\in{\cal C}^{1,3}_K([0,\tf ]\times\R^d).$ $\triangle$ stands the Laplacian operator.

\vk The following notation is devoted to the definition of the rate
functional of the large deviation principle for $(\pi_.^N)_{0\leq
t\leq\tf }$.

\vk For $H\in\cl^{1,2}_K([0,\tf ]\times\R^d)$, we define $\jj\
:\dl([0,\tf ],\mm)\ra\R\cup\{\ty\}$ by
\[
\jj(\pi)=\jj^1(\pi)-\jj^2(\pi)
\]
where
\[
\jj^1(\pi)=\ab u_\tf ,H_\tf \ba-\ab u_0,H_0\ba-\int_0^\tf
\ab u_t,\partial_tH_t\ba\di t,
\]
\[
\jj^2(\pi)=\f{\si}{2}\int_0^\tf
\ab\Phi(u_t),\sum_{i=1}^{d}\bg(\partial_{x_i}^2H_t+(\partial_{x_i}H_t)^2\bg)\ba\di t,
\]

such that $\jj(\cdot)=\ty$ outside $\dl([0,\tf ],\mm)$ or if $\pi_t$ is not
absolutely continuous with respect to the Lebesgue measure $\lm$
for some $0\leq t\leq\tf $.

We are now ready to define the part of the large deviations rate
function, $\ii_0(\cdot): \dl([0,\tf ],\mm)\ra[0,\ty]$ coming from
the stochastic evolution:
\[
\ii_0(\pi)=\sup_{H\in\cl^{1,2}_K([0,\tf ]\times\R^d)}\jj(\pi).
\]

The other part of the large deviations rate function coincides with
the behaviour of deviations coming from the initial state. Let
$h(\cdot|\rho)$ be the entropy defined for a positive function
$\ma:\R^d\ra\R^+$ by
\[
h(\ma|\rho)=\int_{\R^d}\Bg\{
\ma(x)\log\bg(\f{\Phi(\ma(x))}{\Phi(\rho)}\bg)-
\E_m\bg[\log\bg(\f{Z(\Phi(\ma(x))p_0^{-1})}{Z(\Phi(\rho)p_0^{-1})}\bg)\bg]\Bg\}\
\di x.
\]

Thus, the rate function of the large deviation principle is
defined for a density $\rho>0$ by
\[
\ii_\rho(\pi)=\ii_0(\pi)+h(u_0|\rho).
\]

\vk From now on, for each $x\in\Z^d$, we denote by $\eta^l(x)$ the
mean density of particles in a box of length $(2l+1)$ centered at
$x$ : $$\eta^l(x)=\frac{1}{(2l+1)^d}\sum_{|y-x|\leq l}\eta(y).$$

For each cylinder function $\Psi : \X_d\ra\R$, we define
 \Bq \label{tilde}
\tilde{\Psi}(\rho):=m\bg[\nu_{\Phi(\rho)}^p(\Psi)\bg]\ , \Eq
and we say that $\Psi$ is a Lipschitz function if
\[
\exists k_0\in\N\mbox{ and } c_0>0\mbox{ such that }\ \ \
\bg|\Psi(\eta)-\Psi(\xi)\bg|\leq c_0\sum_{|x|\leq
k_0}\bg|\eta(x)-\xi(x)\bg|,
\] 
for all $\eta$ and $\xi$ in $\X_d$.\vk\nit Denote by $\tau_x$ the
shift operator defined by 
$\tau_x\Psi(\eta(\cdot)):= \Psi(\tau_x\eta(\cdot))$ where
$\tau_x\eta(y)=\eta(x+y)$.


We can now state our results:

\newtheorem{theorem}{Theorem}[section]
\begin{theorem} Let $\Psi$ be a cylinder Lipschitz function and
$H\in\cl^{0,2}_K([0,\tf ]\times\R^d)$. Under hypotheses [H1] to
[H4], for all $\dt>0$ we have 
 \bl{th2}
 \Bq
 \bl{superexp}
\lspe\lspN\frac{1}{N^d}\log \pro\Bigg[\Bigg|\int_0^\tf
W_{N,\e}^{H,\Psi}(t,\eta_t)\di t\Bigg|>\dt\Bigg]=-\infty
 \Eq
$m$-almost surely, where

$$W_{N,\e}^{H,\Psi}(t,\eta)=\frac{1}{N^d}\sx
H(t,x/N)\bigg[\tau_x\Psi(\eta)-\tilde\Psi(\eta^{\e N}(x))\bigg].$$

\end{theorem}

This theorem, called the super-exponential estimate, will be a crucial
argument in the proof of the following large deviations principle:


\begin{theorem}
 \bl{th1}
Under hypotheses [H1] to [H4], for every closed subset ${\cal C}$
and every open subset ${\cal O}$ of $\dl([0,\tf ],\mm)$, we have

$$\limsup_{N\ra\ty}\frac{1}{N^d}\log Q^N_{\rho,p}({\cal
C})\leq-\inf_{\pi\in{\cal C}}\ii_\rho(\pi)$$ and
$$\lifN\frac{1}{N^d}\log Q^N_{\rho,p}({\cal O})\geq-\inf_{\pi\in{\cal
O}\cap{\cal A}}\ii_\rho(\pi)$$ $m$-almost surely.
\end{theorem}

\nit{\bf\Large Remarks }

Before starting to prove our results, we would like to mention
some facts and claims that will be used and whose proofs are
omitted. For more details the reader is refered to
Kipnis-Landim's book (1999) and Benois \& al. (1995).

{\bf [R1]\ \ \ } From Lemma I.3.5 of Kipnis-Landim's book (1999),
the function defined by $\vp\longrightarrow\nu_\vp$ for $\vp>0$,
is an increasing function (see also the proof of lemma 4.3 in
Benois \& al. (1995)). Therefore, assumption [H4] implies that for
a fixed environment $p$ defined in the beginning of the last
section, for all $x\in\Z^d$ and $\vp>0$, there exists $\ta:=\ta(x,\vp)>0$ such
that
\[
\nu_\vp^p\bigg[\exp\big\{\ta\om(\eta(x))\big\}\bigg]<\infty\quad m
\mbox{-almost surely.}
\]

{\bf [R2]\ \ \ } Assumption [H4] ensures that the function $\om^*$
defined by (\ref{ome}) is also a continuous convex function such
that $\om^*(0)=0$. \par

{\bf [R3]\ \ \ } A simple computation shows that from the second
condition in [H4], for every $\e>0$ the function $\om^{-1}(r)-\e
r$ is negative for each $r\geq C_2(\e)$, for some constant
$C_2(\e)$ dependent only on $\e$.\par

{\bf [R4]\ \ \ } By definition of $\om$ in [H4], the function
defined on $\R^*_+$ by $\Omega(r)=\f{\om(r)}{r}$ is an increasing
function.\par

{\bf [R5]\ \ \ } For each cylinder Lipschitz function
$\Psi(\cdot)$, the function $\tilde\Psi(\cdot)$ given by
(\ref{tilde}) is also a Lipschitz function (see Lemma I.3.6 of
Kipnis-Landim (1999)). Moreover one can check that $\tilde\Psi(k)\leq Ck$ for all $k\in\Z$ for some constant
$C$.\par

The strategy we adopted to prove the results is similar to the one
presented in Benois \& al. (1995). However, we need some arguments
developed in Koukkous (1999) in order to overcome the lack of
translation invariance of the invariant measures for the zero
range process in random media. We will thus focus only on the main
differences.
\vk From now on, to keep the notation simple, we will restrict our
study to the one-dimensional case. The reader can extend the proofs to
any dimension without any difficulty.\par
\section{\bf Proof of Theorem \ref{th2}}
 \bl{pvsupexpo}
Let $G$ be a positive continuous function on $\R$ defined by \Bq
 \bl{GG}
G(x)=\sup_{y\in[x-1,x+1]}\max\bg\{|H(y)|,|\partial_yH(y)|,|\partial^2_yH(y)|\bg\}.
\Eq 
We have
\begin{eqnarray}
\bl{pr1} \pro\bigg[\int_0^\tf W_{N,\e}^{H,\Psi}(t,\eta_t)\di
t>\dt\Bigg]\nonumber\\ &\ &\hskip-3cm\leq\pro\bigg[\int_0^\tf
\bigg\{W_{N,\e}^{H,\Psi}(t,\eta_t)\di t-\f{\bt}{N}\sx
G\bigg(\f{x}{N}\bigg)\om(\eta_t(x))\bigg\}\di
t>\dt/2\bigg]\nonumber\\ &\ &\hskip0.9cm+\pro\bigg[\int_0^\tf
\f{\bt}{N}\sx G\bigg(\f{x}{N}\bigg)\om(\eta_t(x)) \di
t>\dt/2\bigg]
\end{eqnarray}
for every $\bt>0$.\\
By Tchebycheff exponential inequality the first term in the left
hand side in (\ref{pr1}) is bounded above by
$$\exp\{-N\ta\dt/2\}\esp\Bigg[\exp\ta\int_0^\tf
\bigg\{N\Ww(t,\eta_t)-\bt\sx
G\bigg(\f{x}{N}\bigg)\om(\eta_t(x))\bigg\}\di t\Bigg]$$ for every
$\ta>0$.\\ Therefore, we have to prove two Lemmas:
\newtheorem{lemme}{Lemma}[section]
\begin{lemme}
 \bl{L1}
For every $G\in{\cal C}_K(\R)$,
 \Bq
 \bl{l1}
\lspA\lspN\fN\log\pro\bigg[\int_0^\tf \fN\sx G(x/N)\om(\eta_t(x))
\di t>A\bigg]=-\ty
 \Eq
$m$-almost surely.
\end{lemme}

\begin{lemme}
 \bl{L2}
For any $\ta>0$ and $\bt>0$
 \Bq
 \bl{l2}
\lspe\lspN\fN\log\esp\Bigg[\exp\ta\int_0^\tf
\bigg\{\Ww(t,\eta_t)-\bt\sx G\bigg(\f{x}{N}\bigg)\om(\eta_t(x))\bigg\}\di
t\Bigg]=0.
 \Eq
$m$-almost surely.
\end{lemme}
{\bf Proof of Lemma \ref{L1}}.\\ Using respectively Tchebycheff
exponential inequality and Jensen inequality, we show that for
every positive constant $\ta$, the logarithmic term in (\ref{l1})
is bounded above by
 $$-\ta AN+\log\esp\Bg[\f{1}\tf \int_0^\tf
\exp\bigg\{\sx \ta \tf G(x/N)\om(\eta_t(x))\bigg\}\di t\Bg].$$
From the begining of [R1] and since the product measure $\bar\nu_\rho^p$ is invariant for the
process and $p_x\in[a_0,a_1]$, a simple computation shows that the
right hand side term in (\ref{l1}) is bounded above by
 \Bq
 \bl{ep2}
\lspA\lspN\inf_{\ta>0}\Bigg\{-\ta A+\fN\sx
\log\nu_{\Phi(\rho)a_0^{-1}}\bigg[\exp\bigg\{\ta\tf
G(x/N)\om(\eta(0))\bigg\}\bigg]\Bg\}.
 \Eq
Let $B>0$ be such that $$\mbox{supp}G\subset[-B,B].$$
From {\bf [H4]}, there exists $\ta_0>0$ such that
$$\nu_{\Phi(\rho)a_0^{-1}}\Bg[\exp\bigg\{\ta_0
\tf\|G\|_\infty\om(\eta(0))\bigg\}\Bg]<\ty.$$ The lemma is proved
in fact that (\ref{ep2}) is bounded above by
$$\lspA\bigg\{-\ta_0A+2B\log\nu_{\Phi(\rho)a_0^{-1}}\bigg[e^{\big\{\ta_0\tf\|G\|_\ty\om(\eta(0))\big\}}\bigg]\bigg\}.\hspace{3cm}$$

\nit {\bf Proof of Lemma \ref{L2}}.\\ Let
$$V(\eta)=\ta\bigg\{N\Ww(0,\eta)-\bt\sx
G\bigg(\f{x}{N}\bigg)\om(\eta(x))\bigg\}.$$ Let ${\cal L}^p_V$ be
the generator $N^2{\cal L}_p+V$ and  ${\cal L}^{p,*}_V$ its
adjoint operator, which is equal to $N^2{\cal L}^*_p+V$. If we
denote by $S_{t}^{V,p}$ the semigroup associated to the generator
${\cal L}^p_V$, by the Feyman-Kac formula the expectation in the
lemma is equal to $$\big<S_\tf ^{V,p}1,1\big>\leq\big<S_\tf
^{V,p}1,S_\tf ^{V,p}1\big>^\f{1}{2}.$$ Now, if we denote by $\ld_V$
the largest eigenvalue of the self-adjoint operator ${\cal
L}^p_V+{\cal L}^{p,*}_V$, $$\partial_t\big<S_t ^{V,p}1,S_t
^{V,p}1\big> = \big<({\cal L}^p_V+{\cal L}^{p,*}_V)S_t ^{V,p}1,S_t
^{V,p}1\big>\leq\ld_V\big<S_t ^{V,p}1,S_t ^{V,p}1\big>.$$ By
Gronwall's lemma we show that
 \Bq
 \bl{vpr}
\big<S_\tf ^{V,p}1,S_\tf ^{V,p}1\big>\leq\exp\bigg\{\tf
\ld_V\bigg\}.
 \Eq
Recall that we did not assume $T(\cdot)$ to be symmetric and
therefore  $\nu_{\Phi(\rho)}^p$ can be non-reversible for the
process. However, at this level, our study is dealing with the
reversible generator $N^2({\cal L}_p+{\cal L}_p^*)$. Thus we can
assume the generator ${\cal L}_p$ to be reversible and $T(\cdot)$ given by $T(x)=(1/2)\1_{\{|x|=1\}}$.\par\nit 
Let
$$I_{x,x+1}^p(f)=\f{1}{2}\int
p_xg(\eta(x))\left[\sqrt{f(\eta^{x,x+1})}-\sqrt{f(\eta)}\right]^2\nbr^p(\di\eta),$$

and $D_p(\cdot)$ the Dirichlet form given by $$D_p(f)=\sx
I_{x,x+1}^p(f).$$
Using the variational formula for the largest eigenvalue
of a self-adjoint operator (see appendix A3.1 of Kipnis-Landim (1999)), from (\ref{vpr}) we reduce the proof of the lemma to show
that for every positive $\ta$

$$\lspe\lspN\sup_f\Bigg\{\int\ta\bigg[\Ww(\eta)- \f{\bt}{N}\sx
G\bigg(\f{x}{N}\bigg)\om(\eta(x))\bigg]f(\eta)\nbr^p(\di\eta)-ND_p(f)\Bigg\}\leq0.$$

The supremum is taken over all positive densities functions with
respect to $\nbr^p$. \par\nit We use now some computations from Benois \& al. (1995) and Kipnis \& al. (1989). Let
$$W_l^\Psi(\eta)=\fl\sum_{|y|\leq
l}\tau_y\Psi(\eta)-\tilde\Psi(\eta^l(0))$$ In this way, we can
rewrite the term $$\Ww(\eta)- \f{\bt}{N}\sx
G\bigg(\f{x}{N}\bigg)\om(\eta(x))$$ as
$$\hskip-3cm\fN{\sum_x}\Bigg\{H\bg(\f{x}{N}\bg)\Bg[\tau_x\Psi(\eta)-\fl\sum_{|y-x|\leq
l}\tau_y\Psi(\eta)\Bg]-
\f{\bt}{3}G\bg(\f{x}{N}\bg)\om(\eta(x))\Bg\}$$ $$\hskip-5cm
+\fN{\sum_x}\Bigg\{H\bg(\f{x}{N}\bg)\tau_xW_l^{\Psi}(\eta)-
\f{\bt}{3}G\bg(\f{x}{N}\bg)\om(\eta(x))\Bg\}$$ $$
+\fN{\sum_x}\Bigg\{H\bg(\f{x}{N}\bg)\Bg[\tilde\Psi(\eta^l(x))-\tilde\Psi(\eta^{\e
N}(x))\Bg]- \f{\bt}{3}G\bg(\f{x}{N}\bg)\om(\eta(x))\Bg\}.$$
From the assumption on $\Psi$, we chek easily that there exist $C(\Psi,p)$ such that for all $x\in\Z$ $\Psi(\eta(x))\leq C(\Psi,p)\eta(x)$. Then from the definitions of $\om^*(\cdot)$ and
$G(\cdot)$ (cf. (\ref{ome}) and (\ref{GG})), the first term in the
last expression is bounded above by
$$\hskip-2cm\fN{\sum_x}\Bigg\{\bg|\fl\sum_{|y-x|\leq
l}H\bg(\f{y}{N}\bg)-H\bg(\f{x}{N}\bg)\bg|\Psi(\eta(x))-
\f{\bt}{3}G\bg(\f{x}{N}\bg)\om(\eta(x))\Bg\}$$
$$\hskip-5cm\leq\f{\bt}{3N}\sx G\fxn\Bg\{\f{3{\cal
C}(\Psi,p)l}{\bt N}\eta(x)-\om(\eta(x))\Bg\}$$
$$\hskip-6cm\leq\om^*\Bg\{\f{3{\cal C}(\Psi,p)l}{\bt N}\Bg\}\f{\bt\|G\|_\ty}{3}.$$ 
This last term vanishes as
$N\uparrow\ty$ since $\om^*(\cdot)$ is continuous and
$\om^*(0)=0$.

Now, to achieve the proof of the lemma \ref{L2}, we shall prove:
\begin{lemme}
\label{L3} 
For any $b>0$
$$\hskip-13cm\lspl\lspN\sup_f$$ \Bq \label{l3}
\Bigg\{\fN{\sum_x}\int\bigg[H\bg(\f{x}{N}\bg)\tau_xW_l^{\Psi}(\eta)-
\bt
G\bg(\f{x}{N}\bg)\om(\eta(x))\bg]f(\eta)\di\nbr^p(\di\eta)-bND_p(f)\Bigg\}\leq0\Eq
m-almost surely. The supremum is taken over all positive densities
functions with respect to $\nbr^p$.
\end{lemme}
And, thanks to remarks [R5], we have to prove that:
\begin{lemme}
\label{L4} 
For any $b>0$
$$\hskip-11cm\lspl\lspe\lspN\sup_f$$ \Bq \label{l4}
\Bigg\{\fN{\sum_x}\int\bigg[H\bg(\f{x}{N}\bg)\bg|\eta^{\e
N}(x)-\eta^{l}(x)\bg|-\bt
G\bg(\f{x}{N}\bg)\om(\eta(x))\bg]f(\eta)\di\nbr^p(\di\eta)-bND_p(f)\Bigg\}\leq0\Eq
m-almost surely. The supremum is taken over all positive densities
functions with respect to $\nbr^p$.
\end{lemme}
{\bf Proof of Lemma \ref{L3}}.\\ Using the convexity of $\om$ and
definition of $G$, we check that
 \begin{eqnarray}
\label{eq33.1}
\fN\sx\bg|H\fxn\bg|\om(\eta^l(x))&\leq&\fN\sx\bg|H\fxn\bg|\fl\sum_{|y-x|\leq
l}\om(\eta(y))\nonumber\\ &=&\fN\sx\om(\eta(x))\fl\sum_{|y-x|\leq
l}\bg|H(y/N)\bg|\nonumber\\ &\leq& \fN\sx\om(\eta(x))G\fxn
\end{eqnarray}
At the beginning, we introduce some notations in order to deal in our
study of (\ref{l3}) with the boxes of length $(2l+1)$. Indeed, the
term
$$H\bg(\f{x}{N}\bg)\tau_xW_l^{\Psi}(\eta)-\bt\bg|H\fxn\bg|\om(\eta^l(x))$$
depends on $\eta$ only through $\eta(x-l)\cdots\eta(x+l)$. Thus we
may restrict the integral to microscopic blocks.
Denote by $\wg_l=\{-l\cdots\l\}$ the box of length $(2l+1)$ centered at
the origin. For a fixed $z\in\Z$, we denote by $\wg_{z,l}$ the box
$z+\wg_l$, by $\X^l$ the configuration space $\N^{\wg_l}$,
 by $\bar\nu_{\rho,z,l}^p$ the product measure $\nbr^{\ta_zp}$
restricted to $\X^l$, by $f_{z,l}$ the density, with respect to
$\bar\nu_{\rho,z,l}^{p}$, of the marginal of the measure
$f(\eta)\nbr^{\ta_zp}(d\eta)$ on $\X^l$ and by $D_{\rho,z,l}^p(h)$
the Dirichlet form on $\X^l$ given by
$$D_{\rho,z,l}^p(h)=\sum_{|x-y|=1\atop x,y \in \wg_{z,l}}\int
p_xg(\eta(x))\left[\sqrt{h(\eta^{x,y})}-\sqrt{h(\eta)}\right]^2\bar\nu_{\rho,z,l}^p(\di\eta).$$

Thus, from (\ref{eq33.1}) and since the Dirichlet form is convex
(by Schwarz inequality), the supremum in the lemma is bounded
above by the supremum over all positive densities $f$ (with
respect to $\nbr^p$) of the term \Bq \label{eq33.2}
\fN{\sum_x}\Bg\{\int\bigg[H\bg(\f{x}{N}\bg)W_l^{\Psi}(\eta)-
\bt\bg|H\fxn\bg|\om(\eta^l(0))\bg]f_{x,l}\bar\nu_{\rho,x,l}^p(\di\eta)-\f{bN^2}{C(l)}D_{\rho,x,l}^p(f_{x,l})\Bg\}
\Eq As in the proof of lemma 3.1 of Koukkous (1999) we may now
characterize the sites $x$ where the environment degenerates
(behaves badly).\\ Fix $\dt>0$ , $\fa>0$ and $n\in\N$ sufficiently
large such that $\frac{a_1-a_0}{n}<\dt$. For $0\leq j\leq n-2$, let
$I_j^\dt=[\bt_j,\bt_{j+1}[$ where $\bt_j\in[a_0,a_1]$  is such
that $$\bt_j = a_0+(a_1-a_0)\Bigg(\frac{j}{n}\Bigg) .$$ Let
$I_{n-1}^\dt=[\bt_{n-1},a_1]$ and notice that, for $0\leq j\leq
n-1$, we have $|\bt_{j+1}-\bt_j|<\dt$. \vskip0.2truecm Fix $k<l$
and $L=[\frac{2l+1}{2k+1}]$. We now subdivide $\wg_l$ into $L$
disjoint cubes of length $(2k+1)$; let $B_1,\cdots,B_L$ be such
that $$B_i\subseteq\wg_l,\hspace{1cm} B_i\cap B_j=\emptyset\
\mbox{\ for\ } i\neq j\mbox{\ \ and}\hspace{1cm}
B_i=x_i+\wg_k\mbox{\ \ for some}\ x_i\in\Z .$$ We take $B_1=\wg_k$
and let $B_0=\wg_l-\cup_{j=1}^LB_j$. Finally we define
$B_j(x)=x+B_j$  for $0\leq j\leq L$\ and\  $x\in\Z$ . \vk For
$x\in\Z$,\ $n\in\N$,\ $0\leq j\leq n-1$ and $1\leq i\leq L$,
$N_{x,j,i}^{l,k,\dt}(p)$ is the average number of sites $y$ in
$B_i(x)$ such that $p_y\in I_j^\dt$:
$$N_{x,j,i}^{l,k,\dt}(p)=\frac{1}{(2k+1)}\sum_{z\in
B_i(x)}\1_{\{p_z\in I_j^\dt\}}.$$ For $\fa>0$, we let
$$\A_{x,i,\fa}=\bigg\{p,\hspace{4mm}
\bigg|N_{x,j,i}^{l,k,\dt}(p)-m(I_j^\dt)\bigg|\leq\fa\ \  \mbox{for
all}\ j,\ 0\leq j\leq n-1\bigg\}.$$ To keep notation simple, we
denote $\A_{0,1,\fa}$ by $\A_\fa$. Let $$\A_{x,\fa}=\bigg\{p ,
\hspace{3mm} \frac{1}{L}\sum_{i=1}^L\1_{\{p\in \A_{x,i,\fa}\}}\geq
1-\fa\bigg\}.$$
From the definition of $\om^*$ and the property of $\Psi(\cdot)$ and
$\tilde\Psi(\cdot)$ given in the remarks [R5], a simple
computation shows that the integral term in (\ref{eq33.2}) is
bounded by \Bq \label{eq33.3}
C_1=\bt\|H\|_\ty\om^*\bg(\f{2C(\Psi,p)}{\bt}\bg). \Eq

Therefore, the supremum over all positive densities $f$ (with
respect to $\nbr^p$)  of the term (\ref{eq33.2}) is bounded above by
$$\f{1}{N}\sx\sup_{p\in\A_{0,\fa}}\sup_{h\in\Cf_p^l}\Bg\{\int\bigg[H\bg(\f{x}{N}\bg)W_l^{\Psi}(\eta)-\bt\bg|H\fxn\bg|\om(\eta^l(0))\bg]h(\eta)\nbrl
(\di\eta)-\f{bN^2}{C(l)}D_{\rho,0,l}^p(h)\Bg\}$$ 
\Bq
\label{eq33.4} \hskip-9cm+\
C_1\f{1}{N}\sx\1_{\{p\notin\A_{x,\fa}\}} 
\Eq 
where $\Cf_p^l$ is the set of positive density functions with respect to $\nbrl$.\\ 
By ergodicity and stationary of the environment law, the
second term converges $m$-almost surely, as $N\uparrow\ty$, to
$$C_1m\bg\{p\notin\A_{0,\fa}\bg\}.$$ Again the ergodicity of $m$
ensures that this expression vanishes as $l\uparrow\ty$ and
$k\uparrow\ty$ afterwards. Now, let us turn to the first term in
(\ref{eq33.4}). If we denote 
$$\E_h^p[f]=\int h(\eta)f(\eta)\di
\nbrl (\eta),$$ 
the integral term in (\ref{eq33.4}) is bounded
above by
$$2C(\Psi)\bg|H\fxn\bg|\Bg\{\E_h^p\bg[\eta^l(0)\bg]-\f{\bt}{2C(\Psi)}\E_h^p\bg[\om(\eta^l(0))\bg]\Bg\}.$$

Recall that $\om$ is a convex and increasing function. Thus, by
Jensen's inequality, the last expression is bounded above by
$$2C(\Psi)\bg|H\fxn\bg|\Bg\{\om^{-1}\Bg[\E_h^p\bg[\om(\eta^l(0))\bg]\Bg]-\f{\bt}{2C(\Psi)}\E_h^p\bg[\om(\eta^l(0))\bg]\Bg\}.$$

From the remarks  [R3], we claim that there exists a finite constant
$C_2=C_2(\bt,C(\Psi))$ such that the integral term in
(\ref{eq33.4}) is negative if $\E_h^p\bg[\eta^l(0)\bg]\geq C_2$.
\nit
Let $B>0$ be such that $suppH\subset[-B,B]$, then from (\ref{eq33.3})
and the last claim, we check that the first term in (\ref{eq33.4})
is bounded above by 
$$(2B+1)\|H\|_\ty\sup_{p\in
\A_{0,\fa}}\sup_{f\in{\Cf_l^p(\f{2C_1C(l)}{bN^2},C_2)}}\Bg|\int
W_l^{\Psi}(\eta)f(\eta)\nbrl (\di\eta)\Bg|$$ 
where $\Cf_l^p(a,b)$
is defined for positive constant $a$ and $b$ by
$$\Cf_l^p(a,b)=\Bg\{f\in\Cf_l^p : D_{\rho,0,l}^p(f)\leq a \
\mbox{and}\ \  \E_f^p\bg[\om(\eta^l(0))\bg]\leq b\Bg\}.$$ 
The weak
topology of the set of probability measures on $\X^l$ ensures that,
by definition, $\Cf_l^p(\f{2C_1C(l)}{bN^2},C_2)$ is one of its
compact subsets. Therefore, by the lower semi-continuity of the Dirichlet
form, we know that

$$\hskip-2cm\lspN\sup_{p\in
\A_{0,\fa}}\sup_{f\in{\Cf_l^p(\f{2C_1C(l)}{bN^2},C_2)}}\Bg|\int
W_l^{\Psi}(\eta)f(\eta)\nbrl (\di\eta)\Bg|$$ \Bq \label{eq33.5}
\leq\sup_{p\in \A_{0,\fa}}\sup_{f\in{\Cf_l^p(0,C_2)}}\Bg|\int
W_l^{\Psi}(\eta)f(\eta)\nbrl (\di\eta)\Bg|. 
\Eq
From the assumption on $\Psi$ (and $\tilde\Psi$), for every positive
constant $C_3$, the term in absolute value is bounded above by

$$2C(\Psi)\int\1_{\{\eta^l(0)\geq  C_3\}}\eta^l(0)f(\eta)\nbrl
(\di\eta)+\Bg|\int W_l^{\Psi}(\eta)\1_{\{\eta^l(0)\leq
C_3\}}f(\eta)\nbrl (\di\eta)\Bg|.$$

By remarks [R4], the first term in the last expression is bounded
above by
\begin{eqnarray}
2C(\Psi)\bg(\f{C_3}{\om(C_3)}\bg)\int\om(\eta^l(0))f(\eta)\nbrl(\di\eta)&=&2C(\Psi)\bg(\f{C_3}{\om(C_3)}\bg)\E_f^p\bg[\om(\eta^l(0))\bg]\nonumber\\
&\leq&2C_2C(\Psi)\bg(\f{C_3}{\om(C_3)}\bg)\nonumber
\end{eqnarray}
for all $f\in\Cf_l^p(0,C_2)$. From (H4), this last term vanishes
as $C_3\uparrow\ty$. At this point, we achieve by proving that \Bq \label{eq33.6} \lspk\lspl\sup_{p\in
\A_{0,\fa}}\sup_{f\in{\Cf_l^p(0,C_2)}}\Bg|\int
W_l^{\Psi}(\eta)\1_{\{\eta^l(0)\leq C_3\}}f(\eta)\nbrl
(\di\eta)\Bg|\leq {\cal C}(\dt,\fa) 
\Eq 
where ${\cal C}(\dt,\fa)$
vanishes as $\fa\downarrow0$ and $\dt\downarrow 0$ afterwards. We
omit this proof since it is developed in the proof of
lemma 3.1 in Koukkous (1999).\par\nit {\bf Proof of Lemma \ref{L4}}.

First of all, we approximate (replace) the average over a small
macroscopic box by an average over large microscopic boxes. More
precisely, for $N$ sufficiently large we check that
\begin{eqnarray}
\ &\ &\hskip-1.5cm\fN\sx\bg|H\fxn\bg|\bg|\eta^{\e
N}(x)-\eta^{l}(x)\bg|\nonumber\\
&\leq&\fN\sx\bg|H\fxn\bg||\feN\sum_{2l+1<|y|\leq\e
N}|\eta^l(x)-\eta^{l}(x+y)|+{\cal O}\bg(\f{l}{\e N}\bg)\sx
G\fxn\eta(x)\nonumber\\
&\leq&\fN\sx\bg|H\fxn\bg|\feN\sum_{2l+1<|y|\leq\e
N}|\eta^l(x)-\eta^{l}(x+y)|+\f{\bt}{N}\sx
G\fxn\om(\eta(x))\nonumber
\end{eqnarray}
Define
$$\om_l(\eta,\xi,x,z)=\Bg(\om(\eta^l(x))+\om(\xi^l(z))\Bg)$$
$$W_A^l(\eta,\xi,x,z)=|\eta^l(x)-\xi^{l}(y)|\1_{\{\eta^l(x)\vee\xi^l(z)\leq
A\}}$$ and to keep notation simple, we denote
$W_A^l(\eta,\xi,0,0)$ by $W_A^l(\eta,\xi)$ and
$\om_l(\eta,\xi,0,0)$ by $\om_l(\eta,\xi)$.

As in the previous proof, we introduce an indicator function and
in the same way as in (\ref{eq33.1}), we reduce our proof to show
that, for every positive constant $A$ \Bq \label{eq34.1}
\hskip-4cm\lspl\lspe\lspN\sup_f\Bg\{\fN\sx\bg|H\fxn\bg|\feN\sum_{2l+1<|y|\leq\e
N} \Eq
$$\hskip3cm\int\Bg[W_A^l(\eta,\eta,x,x+y)-\bt\om_l(\eta,\eta,x,x+y)\Bg]f(\eta)\nbr(\di\eta)-bND_p(f)\Bg\}\leq
0$$
From the definition of
$$ \feN\sum_{2l+1<|y|\leq\e N}W_A^l(\eta,\eta,x,x+y)$$ and since $\eta^l(x)$ and $\eta^l(x+y)$ depend on the configuration
$\eta$ only through its values on the set
$$\wg_{x,y,l}:=\wg_{x,l}\cup \Big(y+\wg_{x,l}\Big),$$ we shall
replace $f$ by its conditional expectation with respect to the
$\sigma\mbox{-algebra}$ generated by $\{\eta(z);
z\in\wg_{x,y,l}\}$. Some notation are necessary. For all $y\in\Z$,
we define the shift operator $\ta_y(\cdot)$  on environments
by $(\ta_yp)(x)=p(x+y)$.

For fixed integer $l$ and environments $p$ and $q$, we denote by
$\tilde\X^l$ the configuration space $\N^{\wg_l}\times\N^{\wg_l}$,
by $\bar\nu_{\rho,x,l}^{p,q}$
 the product measure $\bar\nu_\rho^{\ta_xp}\
\otimes\bar\nu_\rho^{\ta_xq}$ restricted to $\tilde\X^l$, and by
$\fk$ the conditional expectation of $f$ with respect  to the
$\sigma\mbox{-algebra}$ generated by $\{\eta(z);
z\in\wg_{x,y,l}\}$. Thus the supremum in (\ref{eq34.1}) is bounded
above by
$$\hskip-8cm\sup_f\Bg\{\fN\sx\bg|H\fxn\bg|\feN\sum_{2l+1<|y|\leq\e
N}$$
$$\hskip3cm\int\Bg[\tau_xW_A^l(\xi_1,\xi_2)-\bt\om_l(\xi_1,\xi_2)\Bg]f_{x,y,l}^p(\xi_1,\xi_2)\nbrxl^{p,\ta_yp}(\di\xi)-bND^p(f)\Bg\}.$$

Let us turn now to the Dirichlet form of $f_{x,y,l}^p$ into microscopic
boxes $\wg_{x,y,l}$. Let $D_l^{p,q}(h)$ be
$$D_l^{p,q}(h)=I_{l,1}^{p,q}(h)+I_{l,2}^{p,q}(h)+\sum_{z,z'\in
\wg_l\atop |z-z'|=1}I_{z,z',1}^{p,q}(h)+\sum_{z,z'\in\wg_l\atop
|z-z'|=1}I_{z,z',2}^{p,q}(h)$$ where, for each $ z,z' \in\wg_l$,
such that $|z-z'|=1$, $$I_{z,z',1}^{p,q}(h)=1/2\int
p_zg(\xi_1(z))\Big[\sqrt{h(\xi_1^{z,z'},\xi_2)}
-\sqrt{h(\xi_1,\xi_2)}\Big]^2\nrl^{p,q}(d\xi),$$
$$I_{z,z',2}^{p,q}(h)=1/2\int
q_zg(\xi_2(z))\Big[\sqrt{h(\xi_1,\xi_2^{z,z'})}
-\sqrt{h(\xi_1,\xi_2)}\Big]^2\nrl^{p,q}(d\xi),$$
$$I_{l,1}^{p,q}(h)=1/2\int
p_0g(\xi_1(0))\Big[\sqrt{h(\xi_1^{0,-},\xi_2^{0,+})}
-\sqrt{h(\xi_1,\xi_2)}\Big]^2\nrl^{p,q}(d\xi),$$
$$I_{l,2}^{p,q}(h)=1/2\int
q_0g(\xi_2(0))\Big[\sqrt{h(\xi_1^{0,+},\xi_2^{0,-})}
-\sqrt{h(\xi_1,\xi_2)}\Big]^2\nrl^{p,q}(d\xi).$$ The
configurations $\xi^{0,\pm}(\cdot)$ are defined by
$$\xi^{0,\pm}(z) = \left\{\begin{array}{ll}
                       \xi(z) & \ \mbox{if $z\neq 0$ } \\
                       \xi(0)\pm 1& \ \mbox{if $z=0$.}
                         \end{array}
                \right.$$

We claim that \Bq \label{eq34.2} \fN\sx\feN\sum_{2l+1<|y|\leq\e
N}D_l^{p,\ta_y p}\bg(f_{x,y,l}^p\bg)\leq C(l)\e^2ND^p(f). \Eq The
proof of the claim is omitted. See Lemma 4.3 of Koukkous (1999)
for more details.\\
From the same notation in the proof of Lemma \ref{L3}, we separate the
sites where the environment behaves badly and repeat the computation in the begining of
(\ref{eq33.3}). Using (\ref{eq34.2}) and introducing the
indicator function of the environements afterwards, our lemma is a consequence of the following results
\begin{lemme}
\label{LL3} $$\lspl\lspe\lspN\feN\sum_{2l+1<|y|\leq\e
N}\fN\sx\Bg[\1_{\{\ta_xp\notin\A_{0,\fa}\}}+\1_{\{\ta_{x+y}p\notin\A_{0,\fa}\}}\Bg]=0$$
$m$ almost surely.
\end{lemme}
\begin{lemme}
\label{LL4} 
For positive constants $a$ and $b$, let
$$\Cf^{p,q}_l(a,b)=\Bg\{h\geq0,\E_{{\nrl^{p,q}}}[h]=1,
D_l^{p,q}(h)\leq a, \E_h^{p,q}\bg[\om_l(\xi_1,\xi_2)\bg]\leq
b\Bg\}$$ \Bq \label{eq34.3}
\lspk\lspl\lspe\lspN\sup_{p,q\in\A_{0,\fa}}\sup_{h\in\Cf^{p,q}_l\bg(\f{\eN}{bN^2}C_1,C_2\bg)}\E_h^{p,q}\Bg(W_A^l(\xi_1,\xi_2)\Bg)\leq{\cal
C}(\dt,\fa) \Eq where ${\cal C}(\dt,\fa)$ vanishes as
$\fa\downarrow0$ and $\dt\downarrow 0$ afterwards.
\end{lemme}

The lemma \ref{LL3} is trivially proved using the ergodicity and
stationarity of $m$. (see Koukkous (1999)).\par\nit Since
$\Cf^{p,q}_l\bg(\f{\eN}{bN^2}C_1,C_2\bg)$ is a compact subset of
the probability measures set on $\X^l\times\X^l$ endowed with the weak topology,
by the lower semi-continuity of the  Dirichlet form, to prove
(\ref{eq34.3}) it is enough to prove that

$$\lspdt\lspfa\lspk\lspl\sup_{p,q\in\A_{0,\fa}}\ \
\sup_{h\in\Cf^{p,q}_l(0,C_2)}\E_h^{p,q}\Bg(W_A^l(\xi_1,\xi_2)\Bg)=0.$$
which is proved in Koukkous (1999) ( see
the proof of lemma 4.2 at formula (23)).
\section{\bf Proof of Theorem \ref{th1}}
The proof of lower bound presented in Benois \& al. (1995) is easily adapted
in this case using some computations already developed in the
previous proof of super-exponential estimate and some arguments
presented in the below upper bound's proof. We therefore omit
details for the reader.\par
 \bl{LDP}
Let $H\in\cl^{1,2}_K([0,\tf ]\times\R)$ and $\ma\in{\cal
C}_\rho(\R)$. From Girsanov formula, the Radon-Nikodym derivative
of $\hpro$ with respect to $\pro$ is given by \Bq \label{girsanov}
\exp
N\Bg\{\jj^1(\pi_t^N)+h_\ma^{p,N}(\pi_0^N|\rho)-N\int_0^t\sum_{x,y}p_xg(\eta_s(x))T(y)\bg[e^{\{H(t,\f{x+y}{N})-H(t,\f{x}{N})\}}\
-\ 1\bg]\di s\Bg\} \Eq where  $h_\ma^{p,N}(\cdot|\rho) :{\cal
M}_+\ra\R$ is defined by
$$h_\ma^{p,N}(\mu|\rho)=\ab\mu,\log\Bg(\f{\Phi(\ma(\cdot))}{\Phi(\rho)}\Bg)\ba-\fN\sum_x\log\Bg[\f{Z(\Phi(\ma(x/N))p_x^{-1})}{Z(\Phi(\rho)p_x^{-1})}\bg].$$

{\bf Upper bound :}\\ The proof is dealing only with a fixed
compact subset ${\cal C}$ of $\dl([0,\tf ],\mm)$. To extend this
result to a closed subset, we need exponential
tightness for $Q^N_{\rho,p}$. It is easily obtained thanks to the
proof presented in Benois (1996) (see also Lemma V.1.5 in
Kipnis-Landim (1999)).\par 
For every $q>1$,
$$Q^N_{\rho,p}(\cl)=\E^N_{\rho,p}\Bg[\Bg(\f{\di\pro}{\di\hpro}\Bg)^{1/q}\Bg(\f{\di\hpro}{\di\pro}\Bg)^{1/q}\1_{\{\pi^N\in\cl\}}\Bg].$$

Let $\vartheta_\e$ be the approximation of identity defined by
$(2\e)^{-1}\1_{[-\e,\e]}(x)$ and $\ast$ the classic convolution
product.\par For $0\leq s\leq\tf $, let

$$u_{\e,N}^{p,H}(\eta_s)=\f{\si}{2N}\sum_k\{\ti_x^2H(s,k/N)+[\ti_xH(s,k/N)]^2\}\{p_kg(\eta_s(k))-\Phi(\eta^{\e
N}_s(k))\}$$

and

$$u_{N,H}^p(\eta_s)=\fN\sum_kp_kg(\eta_s(k))\Bg\{\sum_jT(j)N^2\bg[e^{\{H(t,\f{k+j}{N})-H(t,\f{k}{N})\}}\
-\1\bg]\hspace{3cm}$$
$$\hspace{9cm}-\f{\si}{2}\bg\{\ti_x^2H(s,k/N)+(\ti_xH(s,k/N))^2\bg\}\Bg\}$$

From (\ref{girsanov}), a simple computation shows that
$\bg(\di\pro/\di\hpro\bg)$ is bounded above by $$\exp
N\Bg\{-\jj^1(\pi_\tf^N)+\jj^2(\pi^N\ast\vartheta_\e)-h_\ma^{p,N}(\pi_0^N|\rho)+\int_0^\tf\bg\{u_{\e,N}^{p,H}(\eta_s)+u_{N,H}^p(\eta_s)\bg\}\di
s\Bg\}$$ Thus, $\fN\log Q^N_{\rho,p}(\cl)$ is bounded above by \Bq
\label{ldp2} \f{1}{q}\sup_{\pi\in\cl}\Bg\{-\jj^1(\pi_\tf
^N)+\jj^2(\pi^N\ast\vartheta_\e)-h_\ma^{p,N}(\pi_0^N|\rho)\Bg\}\hspace{6cm}
\Eq
$$\hspace{4cm}+\fN\log\E^N_{\rho,p}\Bg[\Bg(\f{\di\hpro}{\di\pro}\Bg)^{1/q}\exp\Bg\{\f{N}{q}\int_0^\tf\bg(u_{\e,N}^{p,H}(\eta_s)+u_{N,H}^p(\eta_s)\bg)\di
s\Bg\}\Bg]$$ 
Let $\bar H$ be a real continuous function with the same
support as $\sup_t|H_t|$, such that it bounds above  $\sup_{0\leq
t\leq\tf}[|\ti_x^2H_t|+(\ti_xH_t)^2+|H_t|]$.\par
 Let $C_0\in\N$
such that $supp H\subset [0,\tf ]\times[-(C_0-1),(C_0+1)]$. Using H\"older's inequality, we show that, for $q'\in\R$ such that
$(1/q)+(1/q')=1$, the second term in (\ref{ldp2}) is bounded above
by $$
\f{1}{3Nq'}\log\E^N_{\rho,p}\Bg[\exp\bg\{\f{3Nq'}{q}\bg(\int_0^\tf
u_{\e,N}^{p,H}(\eta_s)\di s-\int_0^\tf
\f{\fa}{N}\sum_k\bh(\fkn)\om(\eta_s(k))\di
s\bg)\bg\}\Bg]\hspace{2cm}$$ \Bq \label{ldp3} +
\f{1}{3Nq'}\log\E^N_{\rho,p}\Bg[\exp\bg\{\f{3Nq'}{q}\int_0^\tf
u_{N,H}^p(\eta_s)\di s\bg\}\Bg]\hspace{4cm} \Eq $$\hspace{3cm}+
\f{1}{3Nq'}\log\E^N_{\rho,p}\Bg[\exp\bg\{
\f{3Nq'}{q}\bg(\f{\fa}{N}\int_0^\tf
\sum_k\bh(\fkn)\om(\eta_s(k))\di s\bg)\bg\}\Bg]$$

Using similar arguments as in the proof of lemma \ref{L2} ( see
(\ref{ep2})), we check that the last term in (\ref{ldp3}) is
bounded above by
$$R_1(\fa,q,H)=\f{2C_0}{3q}\log\nu_{\Phi(\rho)a_0^{-1}}\bigg[e^{\big\{\f{3\fa
q'\tf }{q}\|{\bar H}\|_\ty\om(\eta(0))\big\}}\bigg]$$ which 
vanishes as $\fa\downarrow0$ for each fixed $q$ and $H$ thanks to
assumption [H4].\par
From assumption [H2], we check that $g(k)\leq g^*k$ for all $k\in\Z$
and therefore $\Phi(\rho)\leq g^*\rho$. Thus, we repeat the same
argument as above, a simple computation shows that the second
term in (\ref{ldp3}) is bounded above by
$$R_2(q,H,N)=\f{2C_0}{3q'}\log\nu_{\Phi(\rho)a_0^{-1}}\bigg[e^{\big\{\f{\bt}{N}\eta(0)\big\}}\bigg]$$

where $\bt=\bt(\tf,g^*,H,a_1,q,\si)$.

\nit For each fixed $q$ and $H$, it is easy to see that
$R_2(q,H,N)$ vanishes as $N\uparrow\ty$.\par Let us turn to the
first term in (\ref{ldp3}) and denote $R_3(\fa,q,H,\e)$ its limit
when $N\uparrow\ty$. A similar computation as in the proof of
the super-exponential estimate ( see lemma \ref{L2} and its proof),
gives that $$\lim_{\e\rightarrow 0}R_3(\fa,q,H,\e)=0$$ for all
$\fa>0$, $q>1$ and smooth function $H$.   \par 
In the other hand
notice that by a simple computation and from the ergodicity and
stationarity of $m$, we prove that $h_\ma^{p,N}(\pi_0^N|\rho)$
converges (uniformly in $\pi\in\cal C$) to $h(\ma|\rho)$ when
$N\uparrow\ty$.\par We therefore proved that
$\overline{\lim}_{N\ra\ty}(1/N)\log Q^N_{\rho,p}(\cl)$ is bounded
above by
$$\inf_{H,\ma,q,\fa,\e}\Bg\{\f{1}{q}\sup_{\pi\in\cl}\bg\{-\jj^1(\pi)+\jj^2(\pi\ast\vartheta_\e)-h(\ma|\rho)\bg\}+R_3(\fa,q,H,\e)+R_1(\fa,q,H)\Bg\}$$

\nit where the infimum is taken over all $H\in\cl^{1,2}_K([0,\tf
]\times\R)$, $\ma\in{\cal C}_\rho(\R)$, $q>1$, $\fa>0$ and
$\e>0$.\par 
At this level, using the continuity of
$\jj^2(\cdot\ast\vartheta_\e)$ for every fixed $H$ and $\e>0$,
the compacity of $\cal C$ and the arguments developed in
(Kipnis \& al. (1989)) to permute the supremum and infimum, we check
that this last expression is bounded above by

$$-\inf_{\pi\in\cl}\sup_{H,\ma,q,\fa,\e}\Bg\{\f{1}{q}\bg\{-\jj^1(\pi)+\jj^2(\pi\ast\vartheta_\e)-h(\ma|\rho)\bg\}+R_3(\fa,q,H,\e)+R_1(\fa,q,H)\Bg\}$$

We conclude therefore our proof by letting $\e\downarrow0$.
$\fa\downarrow0$ and $q\downarrow1$.

\vk\noindent {\bf Acknowledgements.} 
A.K. thanks Fapesp support (FAPESP n. 99/06918-0). This work is part
of FAPESP Tematico n. 95/0790-1, and FINEP Pronex n. 41.96.0923.00.

Instituto de Matem\'atica, Estat\'{\i}stica e Computa\c{c}\~ao
Cient\'{\i}fica

Universidade Estadual de Campinas

Caixa Postal 6065

13083-010 Campinas, SP, Brasil.

{\it koukkous@ime.unicamp.br}

{\it herve@ime.unicamp.br}
\end{document}